\documentclass[twoside,reqno]{amsart}

\usepackage{amsmath}
\usepackage{amssymb}
\usepackage{amsthm}

\newtheorem{theorem}{Theorem}[section]

\newtheorem{proposition}[theorem]{Proposition}

\theoremstyle{definition}

\numberwithin{equation}{section}

\begin{document}

\title[Some identities involving Gegenbauer polynomials]{Some identities involving Gegenbauer polynomials}

\author{Dae San Kim$^1$}
\address{$^1$ Department of Mathematics, Sogang University, Seoul 121-742, Republic of Korea}
\email{dskim@sogang.ac.kr}
\author{Taekyun Kim$^2$}
\address{$^2$ Department of Mathematics, Kwangwoon University, Seoul 139-701, Republic of Korea.}
\email{tkkim@kw.ac.kr}
\author{Seog-Hoon Rim$^3$}
\address{$^3$ Department of Mathematics Education, Kyungpook National University, Taegu 702-701, Republic of Korea.}
\email{shrim@knu.ac.kr}

%\subjclass{}

%\keywords{}

\begin{abstract}
In this paper we derive some interesting identities involving Gegenbauer polynomials arising from the orthogonality of Gegenbauer polynomials for the inner product space ${\mathbb{P}}_n$ with respect to the weighted inner product $<p_1,p_2>=\int_{-1} ^1 p_{1}(x)p_2(x)(1-x^2)^{\lambda-\frac{1}{2}}dx.$
\end{abstract}

\maketitle

\section{Introduction}

The Gegenbauer polynomials are given in terms of the Jacobi polynomials $P_n ^{(\alpha,\beta)}(x)$ with $\alpha=\beta=\lambda-\frac{1}{2}$ $(\lambda>-\frac{1}{2},~\lambda \neq 0)$ by
\begin{equation}\label{1}
\begin{split}
C_n^{(\lambda)}(x)&=\frac{\Gamma\left(\lambda+\frac{1}{2}\right)\Gamma\left(n+2\lambda\right)}{\Gamma(2\lambda)\Gamma\left(n+\lambda+\frac{1}{2}\right)}P_n ^{\left(\lambda-\frac{1}{2},\lambda-\frac{1}{2}\right)}(x) \\
&=\binom{n+2\lambda-1}{n}\sum_{k=0} ^n \frac{\binom{n}{k}(2\lambda+n)_k}{\left(\lambda+\frac{1}{2}\right)_k}\left(\frac{x-1}{2}\right)^k
\end{split}
\end{equation}
where $(a)_k=a(a+1)(a+2)\cdots(a+k-1)$, (see \cite{19, 23}).

From \eqref{1}, we note that $C_k ^{(\lambda)} (x)$ is a polynomial of degree $n$ with real coefficients and $C_n ^{(\lambda)}(1)=\binom{n+2\lambda-1}{n}$. The leading coefficient of $C_n ^{(\lambda)}(x)$ is $2^n \binom{\lambda+n-1}{n}$. By the theory of Jacobi polynomials with $\alpha=\beta=\lambda-\frac{1}{2}$, $\lambda>-\frac{1}{2}$, and $\lambda\neq 0$, we get
\begin{equation}\label{2}
C_n ^{(\lambda)} (-x)=(-1)^n C_n ^{(\lambda)}(x).
\end{equation}
It is not difficult to show that $C_n ^{(\lambda)}(x)$ is a solution of the following Gegenbauer differential equation:
\begin{equation*}
(1-x^2)y^{''}-(2\lambda+1)xy^{'}+n(n+2\lambda)y=0.
\end{equation*}
The Rodrigues' formula for the Gegenbauer polynomials are well known as the following:
\begin{equation}\label{3}
(1-x^2)^{\lambda-\frac{1}{2}}C_n ^{(\lambda)}(x)=\frac{(-2)^n(\lambda)_n}{n!(n+2\lambda)_n}\left(\frac{d}{dx}\right)^n(1-x^2)^{n+\lambda-\frac{1}{2}},~~({\text{see \cite{01, 13}}}).
\end{equation}
The equation \eqref{3} can be easily derived from the properties of Jacobi polynomials.

As is well known, the generating function of Gegenbauer polynomials is given by
\begin{equation}\label{4}
\frac{2^{\lambda-\frac{1}{2}}}{(1-2xt+t^2)^{\frac{1}{2}}(1-xt+\sqrt{1-2xt+t^2})^{\lambda-\frac{1}{2}}}=\sum_{n=0} ^{\infty} \frac{\left(\lambda+\frac{1}{2}\right)_n}{(2\lambda)_n}C_n ^{(\lambda)} (x) t^n.
\end{equation}
The equation \eqref{4} can be also derived from the generating function of Jacobi polynomials.

From \eqref{4}, we note that
\begin{equation}\label{5}
\frac{1}{(1-2xt+t^2)^{\lambda}}=\sum_{n=0} ^{\infty}C_n ^{(\lambda)}(x)t^n,~~(|t|<1,~|x| \leq 1).
\end{equation}
The proof of \eqref{5} is given in the following book: Stein $\&$ Weiss, Introduction to Fourier Analysis in Euclidean space, Princeton University Press, 1971.

By \eqref{1} and \eqref{2}, we get
\begin{equation}\label{6}
\int_{-1} ^1 C_n ^{(\lambda)}(x)C_m ^{(\lambda)}(x)(1-x^2)^{\lambda-\frac{1}{2}}dx=\frac{\pi 2^{1-2\lambda}\Gamma(n+2\lambda)}{n!(n+\lambda)(\Gamma(\lambda))^2}\delta_{m,n}
\end{equation}
where $\delta_{m,n}$ is the Kronecker symbol and it holds for each fixed $\lambda\in {\mathbb{R}}$ with $\lambda>-\frac{1}{2}$ and $\lambda \neq 0$.

The equation \eqref{6} implies  the orthogonality of $C_n ^{(\lambda)}(x)$ and the equation \eqref{6} is important
in deriving our results in this paper. From \eqref{5}, we can derive the following derivative of Gehenbauer polynomials $C_n ^{(\lambda)}(x)$:
\begin{equation}\label{7}
\frac{d}{dx}C_n ^{(\lambda)}(x)=2\lambda C_{n-1} ^{(\lambda+1)}(x),~{\text{for }}n \geq 1.
\end{equation}
By \eqref{7}, we get
\begin{equation}\label{7}
\frac{d^k}{dx^k}C_n ^{(\lambda)}(x)=2^k\lambda^kC_{n-k} ^{(\lambda+k)}(x).
\end{equation}
As is well known, the Bernoulli polynomials $B_n(x)$ are defined by the generating function to be
\begin{equation}\label{9}
\frac{t}{e^t-1}e^{xt}=e^{B(x)t}=\sum_{n=0} ^{\infty}B_n(x)\frac{t^n}{n!},~({\text{see [2-10]}}),
\end{equation}
with the usual convention about replacing $B^n(x)$ by $B_n(x)$. In the special case, $x=0,$ $B_n(0)=B_n$ are called the $n$-th Bernoulli numbers.

From \eqref{9}, we note that
\begin{equation}\label{10}
B_n(x)=(B+x)^n=\sum_{l=0} ^n \binom{n}{l}B_{n-l} x^l,~({\text{see [6-10]}}),
\end{equation}
and
\begin{equation}
B_n ^{'}(x)=\frac{d}{dx}B_n(x)=nB_{n-1}(x).
\end{equation}
The Euler polynomials $E_n(x)$ are also defined by the generating function to be
\begin{equation}\label{11}
\frac{2}{e^t+1}e^{xt}=e^{E(x)t}=\sum_{n=0} ^{\infty}E_n(x)\frac{t^n}{n!},~({\text{see [12-17]}}),
\end{equation}
with the usual convention about replacing $E^n (x)$ by $E_n(x)$. In the special case, $x=0$, $E_n(0)=E_n$ are called the $n$-th Euler numbers. By \eqref{11}, we see that the recurrence formula for $E_n$ is given by
\begin{equation}\label{12}
E_0=1,~(E+1)^n+E_n=2\delta_{0,n},~({\text{see [20-22]}}).
\end{equation}
For each fixed $\lambda \in {\mathbb{R}}$ with $\lambda>-\frac{1}{2}$ and $\lambda \neq 0$, let ${\mathbb{P}}_n=\left\{p(x)\in{\mathbb{R}}[x]~|~\deg p(x) \leq n \right\}$ be inner product space with respect to the inner product
\begin{equation}\label{13}
<p_1(x),p_2(x)>=\int_{-1} ^1 (1-x^2)^{\lambda-\frac{1}{2}}p_1(x)p_2(x)dx,
\end{equation}
where $p_1(x),p_2(x)\in {\mathbb{P}}_n$.

In this paper, we derive some interesting identities involving Gegenbauer polynomials arising from the orthogonality of those for the inner product space ${\mathbb{P}}_n$ with respect to the weighted inner product $<p_1,p_2>=\int_{-1} ^1 p_1(x)p_2(x)(1-x^2)^{\lambda-\frac{1}{2}}dx.$

Our methods used in this paper are useful in finding some new identities and relations on the Bernoulli and  Euler polynomials involving Gegenbauer polynomials.

\section{Some identities involving Gegenabauer polynomials}

Let us take $p(x)=\sum_{k=0} ^nd_kC_k ^{(\lambda)}(x)\in {\mathbb{P}}_n, ~d_k \in {\mathbb{R}}$. Then, by \eqref{6} and \eqref{13}, we get
\begin{equation}\label{14}
\begin{split}
& <p(x),C_k ^{(\lambda)}(x)>=d_k<C_k ^{(\lambda)}(x),C_k ^{(\lambda)}(x)> \\
&=d_k\int_{-1} ^1 (1-x^2)^{\lambda-\frac{1}{2}}C_k ^{(\lambda)}(x)C_k ^{(\lambda)}(x)dx=d_k\frac{\pi2^{1-2\lambda}\Gamma(k+2\lambda ) }{k!(k+\lambda)(\Gamma(\lambda))^2}.
\end{split}
\end{equation}
Thus, from \eqref{14}, we have
\begin{equation}\label{15}
d_k=\frac{(\Gamma(\lambda))^2k!(k+\lambda)}{\pi2^{1-2\lambda}\Gamma(k+2\lambda)}\int_{-1} ^1 (1-x^2)^{\lambda-\frac{1}{2}}p(x)C_k ^{(\lambda)}(x)dx.
\end{equation}
By \eqref{3} and \eqref{15}, we get
\begin{equation}\label{16}
\begin{split}
d_k&=\frac{(\Gamma(\lambda))^2k!(k+\lambda)}{\pi2^{1-2\lambda}\Gamma(k+2\lambda)}\times\frac{(-2)^k(\lambda)_k}{k!(k+2\lambda)_k}\int_{-1} ^1 \left(\frac{d^k}{dx^k}(1-x^2)^{k+\lambda-\frac{1}{2}}\right)p(x)dx \\
&=\frac{(k+\lambda)\Gamma(\lambda)}{(-2)^k\sqrt{\pi}\Gamma(k+\lambda+\frac{1}{2})}\int_{-1} ^1\left(\frac{d^k}{dx^k}(1-x^2)^{k+\lambda-\frac{1}{2}}\right)p(x)dx.
\end{split}
\end{equation}
Therefore, by \eqref{16}, we obtain the following proposition.
\begin{proposition}\label{prop1}
For $p(x) \in {\mathbb{P}}_n$, let
\begin{equation*}
p(x)=\sum_{n=0} ^n d_kC_k ^{(\lambda)}(x),~(d_k \in {\mathbb{R}}).
\end{equation*}Then
\begin{equation*}
d_k=\frac{(k+\lambda)\Gamma(\lambda)}{(-2)^k\sqrt{\pi}\Gamma(k+\lambda+\frac{1}{2})}\int_{-1} ^1\left(\frac{d^k}{dx^k}(1-x^2)^{k+\lambda-\frac{1}{2}}\right)p(x)dx.
\end{equation*}
\end{proposition}
For example, let $p(x)=x^n \in {\mathbb{P}}_n$. From Proposition \ref{prop1}, we note that
\begin{equation}\label{17}
\begin{split}
d_k&=\frac{(k+\lambda)\Gamma(\lambda)}{(-2)^k\sqrt{\pi}\Gamma(k+\lambda+\frac{1}{2})}\int_{-1} ^1\left(\frac{d^k}{dx^k}(1-x^2)^{k+\lambda-\frac{1}{2}}\right)x^ndx \\
&=(-n)\int_{-1} ^1\left(\frac{d^{k-1}}{dx^{k-1}}(1-x^2)^{k+\lambda-\frac{1}{2}}\right)x^{n-1}dx \times \left(\frac{(k+\lambda)\Gamma(\lambda)}{\sqrt{\pi}(-2)^k\Gamma\left(k+\lambda+\frac{1}{2}\right)}\right) \\
&=\cdots \\
&= \frac{(k+\lambda)n!\Gamma(\lambda)}{(n-k)! 2^k \sqrt{\pi} (k+\frac{1}{2}+\lambda)}\int_{-1} ^1 (1-x^2)^{k+\lambda-\frac{1}{2}}x^{n-k}dx\\
&=\left(1+(-1)^{n-k}\right)\frac{(k+\lambda)n!\Gamma(\lambda)}{(n-k)!2^k\sqrt{\pi}\Gamma(k+\frac{1}{2}+\lambda)}\int_{0} ^1(1-x^2)^{k+\lambda-\frac{1}{2}}x^{n-k} dx.
\end{split}
\end{equation}
Let us assume that $n-k\equiv0~({\rm{mod}}~ 2)$. Then, by \eqref{17}, we get
\begin{equation}\label{18}
\begin{split}
d_k&=\frac{(k+\lambda)n!\Gamma(\lambda)}{(n-k)!2^k\sqrt{\pi}\Gamma(k+\frac{1}{2}+\lambda)}B\left(k+\lambda+\frac{1}{2},\frac{n-k+1}{2}\right)\\
&=\frac{\Gamma\left(\frac{n-k+1}{2}\right)\Gamma\left(k+\lambda+\frac{1}{2}\right)}{\Gamma\left(\frac{n+k+2\lambda+2}{2}\right)},
\end{split}
\end{equation}
where $B(\alpha,\beta)$ is the beta function which is defined by $B(\alpha,\beta)=\frac{\Gamma(\alpha)\Gamma(\beta)}{\Gamma(\alpha+\beta)}$.

It is easy to show that
\begin{equation}\label{19}
\begin{split}
\Gamma\left(\frac{n-k+1}{2}\right)&=\frac{n-k-1}{2}\Gamma\left(\frac{n-k-1}{2}\right) \\
& =\left(\frac{n-k-1}{2}\right)\left(\frac{n-k-3}{2}\right)\Gamma\left(\frac{n-k-3}{2}\right)=\cdots\\
&=\frac{\left(\frac{n-k}{2}\right)\left(\frac{n-k-1}{2}\right)\left(\frac{n-k-2}{2}\right)\cdots\frac{2}{2}\Gamma\left(\frac{1}{2}\right)}{\left(\frac{n-k}{2}\right)\left(\frac{n-k-2}{2}\right)\cdots\left(\frac{2}{2}\right)}=\frac{(n-k)!\sqrt{\pi}}{2^{n-k}\left(\frac{n-k}{2}\right)!}.
\end{split}
\end{equation}
Therefore, by \eqref{18} and \eqref{19}, we obtain the following identity:
\begin{equation}\label{20}
x^n=\sum_{0 \leq k \leq n,n-k \equiv0~({\rm{mod}}~2)}\frac{(k+\lambda)n!\Gamma(\lambda)}{2^n\left(\frac{n-k}{2}\right)!\Gamma\left(\frac{n+k+2\lambda+2}{2}\right)}C_k ^{(\lambda)}(x).
\end{equation}
Let us take $p(x)=B_n(x)\in {\mathbb{P}}_n$. Then, by \eqref{10}, we get
\begin{equation}\label{21}
\begin{split}
d_k&=\frac{(k+\lambda)\Gamma(\lambda)}{(-2)^k\sqrt{\pi}\Gamma\left(k+\lambda+\frac{1}{2}\right)}\int_{-1} ^1 \left(\frac{d^k}{dx^k}(1-x^2)^{k+\lambda-\frac{1}{2}}\right)B_n(x)dx \\
&=\frac{(k+\lambda)\Gamma(\lambda)(-n)}{(-2)^k\sqrt{\pi}\Gamma\left(k+\lambda+\frac{1}{2}\right)}\int_{-1} ^1 \left(\frac{d^{k-1}}{dx^{k-1}}(1-x^2)^{k+\lambda-\frac{1}{2}}\right)B_{n-1}(x)dx =\cdots \\
&=\frac{(k+\lambda)\Gamma(\lambda)(-n)(-(n-1))\cdots(-(n-k+1))}{(-2)^k\sqrt{\pi}\Gamma\left(k+\lambda+\frac{1}{2}\right)}\int_{-1} ^1 (1-x^2)^{k+\lambda-\frac{1}{2}}B_{n-k}(x)dx \\
&=\frac{(k+\lambda)\Gamma(\lambda)}{2^k\sqrt{\pi}\Gamma(k+\lambda+\frac{1}{2})}\times\frac{n!}{(n-k)!}\int_{-1} ^1 (1-x^2)^{k+\lambda-\frac{1}{2}}B_{n-k}(x)dx.
\end{split}
\end{equation}
From \eqref{10} and \eqref{20}, we can derive the following equation:
\begin{equation}\label{22}
\begin{split}
& \int_{-1} ^1 (1-x^2)^{k+\lambda-\frac{1}{2}}B_{n-k}(x)dx=\sum_{l=0} ^{n-k}\binom{n-k}{l}B_{n-k-l}\int_{-1}^{1}(1-x^2)^{k+\lambda-\frac{1}{2}}x^l dx \\
&=\sum_{l=0} ^{n-k} \binom{n-k}{l}B_{n-k-l}(1+(-1)^l)\int_0 ^1 (1-x^2)^{k+\lambda-\frac{1}{2}}x^l dx.
\end{split}
\end{equation}
Let us consider that $l\equiv 0~({\rm{mod}}~2)$. Then, by \eqref{22}, we get
\begin{equation}\label{23}
\begin{split}
&\int_{-1} ^1 (1-x^2)^{k+\lambda-\frac{1}{2}}B_{n-k}(x)dx \\
=& 2\sum_{0 \leq l \leq n-k,~l\equiv0~({\rm{mod}}~2)}\binom{n-k}{l}B_{n-k-l}\int_0 ^1 (1-x^2)^{k+\lambda-\frac{1}{2}}x^l dx \\
=& \sum_{0 \leq l \leq n-k,~l\equiv0~({\rm{mod}}~2)}\binom{n-k}{l}B_{n-k-l}\int_0 ^1 (1-y)^{k+\lambda-\frac{1}{2}}y^{\frac{l-1}{2}} dy \\
=& \sum_{0 \leq l \leq n-k,~l\equiv0~({\rm{mod}}~2)}\binom{n-k}{l}B_{n-k-l}\frac{\Gamma\left(k+\lambda+\frac{1}{2}\right)\Gamma\left(\frac{l+1}{2}\right)}{\Gamma\left(\frac{2k+2\lambda+l+2}{2}\right)}.
\end{split}
\end{equation}
For $l \in {\mathbb{Z}}_+$ with $l\equiv0~{\text{mod}}~2)$, we have
\begin{equation}\label{24}
\begin{split}
\Gamma\left(\frac{l+1}{2}\right)&=\Gamma\left(\frac{l-1}{2}+1\right)=\frac{l-1}{2}\Gamma\left(\frac{l-1}{2}\right) \\
=& \left(\frac{l-1}{2}\right)\left(\frac{l-3}{2}\right)\Gamma\left(\frac{l-3}{2}\right)=\cdots \\
=& \left(\frac{l-1}{2}\right)\left(\frac{l-3}{2}\right)...\left(\frac{1}{2}\right)\Gamma\left(\frac{1}{2}\right)
=\frac{\left(\frac{1}{2}\right)^ll!\Gamma\left(\frac{1}{2}\right)}{\left(\frac{l}{2}\right)!}=\frac{l!\sqrt{\pi}}{2^l\left(\frac{l}{2}\right)!}.
\end{split}
\end{equation}
By \eqref{23} and \eqref{24}, we get
\begin{equation}\label{25}
\begin{split}
&\int_{-1} ^1 (1-x^2)^{k+\lambda-\frac{1}{2}}B_{n-k}(x)dx \\
=&\sum_{0 \leq l \leq n-k,~l\equiv 0~({\text{mod}}~2)}\binom{n-k}{l}B_{n-k-l}\frac{\Gamma\left(k+\lambda+\frac{1}{2}\right)\Gamma\left(\frac{l+1}{2}\right)}{\Gamma\left(\frac{2k+2\lambda+l+2}{2}\right)} \\
=& \sum_{0 \leq l \leq n-k,~l\equiv0~({\text{mod}}~2)}\binom{n-k}{l}B_{n-k-l} \frac{l!\sqrt{\pi}}{2^l\left(\frac{l}{2}\right)!}\times\frac{\Gamma\left(k+\lambda+\frac{1}{2}\right)}{\Gamma\left(\frac{2k+2\lambda+l+2}{2}\right)}.
\end{split}
\end{equation}
From \eqref{21} and \eqref{25}, we have
\begin{equation}\label{26}
d_k=\frac{n!(k+\lambda)\Gamma(\lambda)}{2^k(n-k)!}\sum_{0 \leq l \leq n-k,~l\equiv0~({\text{mod}}~2)}\frac{\binom{n-k}{l}B_{n-k-l}l!}{2^l\left(\frac{l}{2}\right)!\Gamma\left(\frac{2k+2\lambda+l+2}{2}\right)}.
\end{equation}
Therefore, by \eqref{26} and Proposition \ref{prop1}, we obtain the following theorem.
\begin{theorem}\label{thm2}
For $n \in {\mathbb{Z}}_+$, we have
\begin{equation*}
\frac{B_n(x)}{n!}=\Gamma(\lambda)\sum_{k=0} ^n \left(\frac{(k+\lambda)}{2^k(n-k)!}\sum_{0 \leq l \leq n-k,~l\equiv0~({\text{mod}}~2)}\frac{\binom{n-k}{l}B_{n-k-l}l!}{2^l\left(\frac{l}{2}\right)!\Gamma\left(\frac{2k+2\lambda+l+2}{2}\right)}\right)C_k ^{(\lambda)}(x).
\end{equation*}
\end{theorem}
By the same method, we get
\begin{equation}\label{27}
\frac{E_n(x)}{n!}=\Gamma(\lambda)\sum_{k=0} ^n \left(\frac{(k+\lambda)}{2^k(n-k)!}\sum_{0 \leq l \leq n-k,~l\equiv0~({\text{mod}}~2)}\frac{\binom{n-k}{l}E_{n-k-l}l!}{2^l\left(\frac{l}{2}\right)!\Gamma\left(\frac{2k+2\lambda+l+2}{2}\right)}\right)C_k ^{(\lambda)}(x).
\end{equation}
From \eqref{1}, we note that
\begin{equation}\label{28}
\begin{split}
&C_{n-k} ^{(\lambda)}(x)C_k ^{(\lambda)}(x) \\
=&\binom{n-k+2\lambda-1}{n-k}\sum_{l=0} ^{n-k}\frac{\binom{n-k}{l}(2\lambda+n-k)_l}{\left(\lambda+\frac{1}{2}\right)_l}\left(\frac{x-1}{2}\right)^l\binom{k+2\lambda-1}{k}\sum_{m-0} ^k \frac{\binom{k}{m}(2\lambda+k)_m}{\left(\lambda+\frac{1}{2}\right)_m}\left(\frac{x-1}{2}\right)^m \\
=&\binom{n-k+2\lambda-1}{n-k}\binom{k+2\lambda-1}{k}\sum_{p=0} ^n\left(\sum_{m=0} ^p \frac{\binom{n-k}{p-m}\binom{k}{m}(2\lambda+k)_m(2\lambda+n-k)_{p-m}}{\left(\lambda+\frac{1}{2}\right)_m\left(\lambda+\frac{1}{2}\right)_{p-m}}\right)\left(\frac{x-1}{2}\right)^p.
\end{split}
\end{equation}
Let us take $p(x)=C_k ^{(\lambda)}(x)C_{n-k} ^{(\lambda)}(x) \in{\mathbb{P}}_n$. From Proposition \ref{prop1}, $p(x)$ can be rewritten as
\begin{equation}\label{29}
p(x)=C_k ^{(\lambda)}(x)C_{n-k} ^{(\lambda)}(x)=\sum_{r=0} ^n d_r C_r ^{(\lambda)}(x),~(d_r \in {\mathbb{R}}).
\end{equation}
Then, by Proposition \ref{prop1} and \eqref{28}, we get
\begin{equation}\label{30}
\begin{split}
d_r=&\frac{(r+\lambda)\Gamma(\lambda)}{(-2)^r\sqrt{\pi}\Gamma\left(r+\lambda+\frac{1}{2}\right)}\int_{-1} ^1 \left(\frac{d^r}{dx^r}(1-x^2)^{r+\lambda-\frac{1}{2}}\right)C_k ^{(\lambda)}(x)C_{n-k} ^{(\lambda)}(x) dx \\
=&\frac{(r+\lambda)\Gamma(\lambda)}{(-2)^r\sqrt{\pi}\Gamma\left(r+\lambda+\frac{1}{2}\right)}\binom{n-k+2\lambda-1}{n-k}\binom{k+2\lambda-1}{k}\\
& \times \sum_{p=0} ^n \left(\sum_{m=0} ^p \frac{\binom{n-k}{p-m}\binom{k}{m}(2\lambda+k)_m}{\left(\lambda+\frac{1}{2}\right)_m\left(\lambda+\frac{1}{2}\right)_{p-m}}(2\lambda+n-k)_{p-m}\right) \\
& \times \int_{-1} ^1 \left(\frac{d^r}{dx^r}(1-x^2)^{r+\lambda-\frac{1}{2}}\right)\left(\frac{x-1}{2}\right)^pdx \\
=& \frac{(r+\lambda)\Gamma(\lambda)}{(-2)^r\sqrt{\pi}\Gamma\left(r+\lambda+\frac{1}{2}\right)}\binom{n-k+2\lambda-1}{n-k}\binom{k+2\lambda-1}{k}\\
& \times \sum_{p=r} ^n \left(\sum_{m=0} ^p \frac{\binom{n-k}{p-m}\binom{k}{m}(2\lambda+k)_m}{\left(\lambda+\frac{1}{2}\right)_m\left(\lambda+\frac{1}{2}\right)_{p-m}}(2\lambda+n-k)_{p-m}\right) \\
& \times \int_{-1} ^1 \left(\frac{d^r}{dx^r}(1-x^2)^{r+\lambda-1}\right)\left(\frac{x-1}{2}\right)^pdx \\
\end{split}
\end{equation}
It is not difficult to show that
\begin{equation}\label{31}
\begin{split}
& \int_{-1} ^1 \left(\frac{d^r}{dx^r}(1-x^2)^{r+\lambda-\frac{1}{2}}\right)\left(\frac{x-1}{2}\right)^pdx
=\frac{(-1)^rp!}{2^p(p-r)!}\int_{-1} ^1 (1-x^2)^{r+\lambda-\frac{1}{2}}(1-x)^{p-r}(-1)^{p-r}dx \\
=& \frac{(-1)^pp!}{2^p(p-r)!}\int_{-1} ^1 (1-x)^{p+\lambda-\frac{1}{2}}(1+x)^{r+\lambda-\frac{1}{2}}dx
= \frac{(-1)^pp!}{2^p(p-r)!}\int_0 ^1 (2-2y)^{p+\lambda-\frac{1}{2}}(2y)^{r+\lambda-\frac{1}{2}}2dy \\
=&\frac{(-1)^p2^{p+\lambda-\frac{1}{2}+r+\lambda-\frac{1}{2}+1}}{2^p}\times\frac{p!}{(p-r)!}\int_0 ^1 (1-y)^{p+\lambda-\frac{1}{2}}y^{r+\lambda-\frac{1}{2}}dy \\
=&(-1)^p2^{r+2\lambda}\frac{p!}{(p-r)!}\times \frac{\Gamma\left(p+\lambda+\frac{1}{2}\right)\Gamma\left(r+\lambda+\frac{1}{2}\right)}{\Gamma\left(r+p+2\lambda+1\right)}.
\end{split}
\end{equation}
From the fundamental theorem of gamma function, we have
\begin{equation}\label{32}
\begin{split}
\frac{\Gamma\left(p+\lambda+\frac{1}{2}\right)}{\Gamma(r+p+2\lambda+1)}&=\frac{\left(p+\lambda-\frac{1}{2}\right)\cdots\left(\lambda+\frac{1}{2}\right)\Gamma\left(\lambda+\frac{1}{2}\right)}{(r+p+2\lambda)\cdots 2\lambda\Gamma(2\lambda)} \\
&=\frac{\left(\lambda+\frac{1}{2}\right)_p\sqrt{\pi}2^{1-\lambda}}{(2\lambda)_{r+p+1}\Gamma(\lambda)}.
\end{split}
\end{equation}
By \eqref{31} and \eqref{32}, we get
\begin{equation}\label{33}
\begin{split}
&\int_{-1} ^1 \left(\frac{d^r}{dx^r}(1-x^2)^{r+\lambda-\frac{1}{2}}\right)\left(\frac{x-1}{2}\right)^pdx \\
=& (-1)^p2^{r+2\lambda}\frac{p!}{(p-r)!}\times\frac{\Gamma\left(p+\lambda+\frac{1}{2}\right)\Gamma\left(r+\lambda+\frac{1}{2}\right)}{\Gamma(r+p+2\lambda+1)} \\
=&(-1)^p2^{r+2\lambda}\frac{p!}{(p-r)!}\Gamma\left(r+\lambda+\frac{1}{2}\right)\times\frac{\left(\lambda+\frac{1}{2}\right)_p\sqrt{\pi}2^{1-\lambda}}{(2\lambda)_{r+p+1}\Gamma(\lambda)}\\
=&(-1)^p2^{r+\lambda +1 }\frac{p!}{(p-r)!}\Gamma\left(r+\lambda+\frac{1}{2}\right)\times\frac{\left(\lambda+\frac{1}{2}\right)_p\sqrt{\pi}}{(2\lambda)_{r+p+1}\Gamma(\lambda)}.\\
\end{split}
\end{equation}
From \eqref{30} and \eqref{33}, we have
\begin{equation}\label{34}
\begin{split}
d_r=&\frac{(r+\lambda)\Gamma(\lambda)}{(-2)^r\sqrt{\pi}\Gamma\left(r+\lambda+\frac{1}{2}\right)}\binom{n-k+2\lambda-1}{n-k}\binom{k+2\lambda-1}{k} \\
&\times\sum_{p=r} ^n\left(\sum_{m=0} ^p \frac{\binom{n-k}{p-m}\binom{k}{m}(2\lambda+k)_m}{\left(\lambda+\frac{1}{2}\right)_m\left(\lambda+\frac{1}{2}\right)_{p-m}}(2\lambda+n-k)_{p-m}\right)\\
& \times (-1)^p2^{r+\lambda+1}\frac{p!}{(p-r)!}\Gamma\left(\lambda+\frac{1}{2}+r\right)\times\frac{\left(\lambda+\frac{1}{2}\right)_p\sqrt{\pi}}{(2\lambda)_{r+p+1}\Gamma(\lambda)} \\
=&(-1)^{r+p}2^{\lambda+1}(r+\lambda)\binom{n-k+2\lambda-1}{n-k}\binom{k+2\lambda-1}{k} \\
&\times \sum_{p=r} ^n\left(\sum_{m=0} ^p \frac{\binom{n-k}{p-m}\binom{k}{m}(2\lambda+k)_m}{\left(\lambda+\frac{1}{2}\right)_m\left(\lambda+\frac{1}{2}\right)_{p-m}}(2\lambda+n-k)_{p-m} \frac{p!\left(\lambda+\frac{1}{2}\right)_p}{(p-r)!(2\lambda)_{r+p+1}}\right).
\end{split}
\end{equation}
Therefore, by \eqref{34}, we obtain the following theorem.
\begin{theorem}\label{thm3}
For $n,k \in {\mathbb{Z}}_+$ with $n \geq k$, we have
\begin{equation*}
\begin{split}
&C_{n-k} ^{(\lambda)}(x)C_k ^{(\lambda)}(x) \\
&=2^{\lambda+1}\binom{n-k+2\lambda-1}{n-k}\binom{k+2\lambda-1}{k}\sum_{r=0} ^n \sum_{p=r}^n\sum_{m=0}^p \left\{(r+\lambda)(-1)^{p+r} \times \right.\\
&\left.\frac{\binom{n-k}{p-m}\binom{k}{m}(2\lambda+k)^m(2\lambda+n-k)_{p-m}p!\left(\lambda+\frac{1}{2}\right)_p}{\left(\lambda+\frac{1}{2}\right)_m\left(\lambda+\frac{1}{2}\right)_{p-m}(p-r)!(2\lambda)_{r+p+1}}\right\}C_r ^{(\lambda)}(x).
\end{split}
\end{equation*}
\end{theorem}

Let us take $p(x)=C_n ^{(\lambda)}(x)\in {\mathbb{P}}_n$. Then, from \eqref{1}, we have
\begin{equation}\label{35}
\begin{split}
C_n ^{(\lambda)}(x)&=\frac{\Gamma\left(\lambda+\frac{1}{2}\right)\Gamma(n+2\lambda)}{\Gamma(2\lambda)\Gamma\left(n+\lambda+\frac{1}{2}\right)}P_n ^{\left(\lambda-\frac{1}{2},\lambda-\frac{1}{2}\right)}(x) \\
&=\frac{(n+2\lambda-1)\cdots(2\lambda)}{\left(n+\lambda-\frac{1}{2}\right)\cdots\left(\lambda+\frac{1}{2}\right)}P_n ^{\left(\lambda-\frac{1}{2},\lambda-\frac{1}{2}\right)}(x)=\frac{\binom{n+2\lambda-1}{n}}{\binom{n+\lambda-\frac{1}{2}}{n}}P_n ^{\left(\lambda-\frac{1}{2},\lambda-\frac{1}{2}\right)}(x).
\end{split}
\end{equation}
In a previous paper, we have shown that
\begin{equation}\label{36}
P_n ^{(\alpha,\beta)}(x)=\sum_{k=0} ^n \binom{n+\alpha}{n-k}\binom{n+\beta}{k}\left(\frac{x-1}{2}\right)^k\left(\frac{x+1}{2}\right)^{n-k},~{\text{(see \cite{11})}}.
\end{equation}
From \eqref{35} and \eqref{36}, we have
\begin{equation}\label{37}
C_n ^{(\lambda)}(x)=\frac{\binom{n+2\lambda-1}{n}}{\binom{n+\lambda-\frac{1}{2}}{n}}\sum_{k=0} ^n\binom{n+\lambda-\frac{1}{2}}{n-k}
\binom{n+\lambda-\frac{1}{2}}{k}
\left(\frac{x-1}{2}\right)^k\left(\frac{x+1}{2}\right)^{n-k},
\end{equation}
and
\begin{equation}\label{38}
\frac{d^k}{dx^k}C_n ^{(\lambda)}(x)=2^k\lambda^kC_{n-k} ^{(\lambda+k)}(x).
\end{equation}
Let $p(x)=C_n ^{(\lambda)} (x)=\sum_{k=0} ^n d_kC_k ^{(\lambda)}(x)$. Then, by Proposition \ref{prop1}, we get
\begin{equation}\label{39}
\begin{split}
d_k&=\frac{(k+\lambda)\Gamma(\lambda)}{(-2)^k\sqrt{\pi}\Gamma\left(k+\lambda+\frac{1}{2}\right)}\int_{-1} ^1 \left(\frac{d^k}{dx^k}(1-x^2)^{k+\lambda-\frac{1}{2}}\right)C_n ^{(\lambda)}(x) dx \\
&=\frac{(k+\lambda)\Gamma(\lambda)}{(-2)^k\sqrt{\pi}\Gamma\left(k+\lambda+\frac{1}{2}\right)}(-1)^k2^k\lambda^k\int_{-1} ^1 (1-x^2)^{k+\lambda-\frac{1}{2}}C_{n-k} ^{(\lambda+k)}(x)dx\\
&=\frac{\lambda^k(k+\lambda)\Gamma(\lambda)}{\sqrt{\pi}\Gamma\left(k+\lambda+\frac{1}{2}\right)}\int_{-1} ^1 (1-x^2)^{k+\lambda-\frac{1}{2}}C_{n-k} ^{(\lambda+k)}(x)dx.
\end{split}
\end{equation}
By \eqref{37}, we get
\begin{equation}\label{40}
\begin{split}
& C_{n-k} ^{(\lambda+k)}(x) \\
=&\frac{\binom{n-k+2(\lambda+k)-1}{n-k}}{\binom{n-k+\lambda+k-\frac{1}{2}}{n-k}}\sum_{l=0} ^{n-k}\binom{n-k+\lambda+k-\frac{1}{2}}{n-k-l}\binom{n-k+\lambda+k-\frac{1}{2}}{l}\left(\frac{x-1}{2}\right)^l\left(\frac{x+1}{2}\right)^{n-k-l} \\
=& \frac{\binom{n+k+2\lambda-1}{n-k}}{\binom{n+\lambda-\frac{1}{2}}{n-k}}\sum_{l=0} ^{n-k}\binom{n+\lambda-\frac{1}{2}}{n-k-l}\binom{n+\lambda-\frac{1}{2}}{l}\left(\frac{x-1}{2}\right)^l\left(\frac{x+1}{2}\right)^{n-k-l}.
\end{split}
\end{equation}
From \eqref{39} and \eqref{40}, we have
\begin{equation}\label{41}
\begin{split}
d_k&=\frac{\lambda^k(k+\lambda)\Gamma(\lambda)}{\sqrt{\pi}\Gamma\left(k+\lambda+\frac{1}{2}\right)}\times\frac{\binom{n+k+2\lambda-1}{n-k}}{\binom{n+\lambda-\frac{1}{2}}{n-k}}\sum_{l=0} ^{n-k}\binom{n+\lambda-\frac{1}{2}}{n-k-l}\binom{n+\lambda-\frac{1}{2}}{l}(-1)^l\left(\frac{1}{2}\right)^{n-k} \\
& \times \int_{-1} ^1 (1-x)^{k+\lambda-\frac{1}{2}+l}(1+x)^{\lambda+n-\frac{1}{2}-l}dx.
\end{split}
\end{equation}
It is easy to show that
\begin{equation}\label{42}
\begin{split}
&\int_{-1} ^1(1-x)^{k+\lambda-\frac{1}{2}+l}(1+x)^{\lambda+n-l-\frac{1}{2}}dx=\int_0 ^1 (2-2y)^{k+\lambda-\frac{1}{2}+l}(2y)^{\lambda+n-l-\frac{1}{2}}2dy \\
=&2^{n+2\lambda+k}\int_0 ^1(1-y)^{k+\lambda+l-\frac{1}{2}}y^{\lambda+n-l-\frac{1}{2}}dy \\
=& 2^{k+n+2\lambda}\frac{\Gamma\left(k+\lambda+l+\frac{1}{2}\right)\Gamma\left(\lambda+n-l+\frac{1}{2}\right)}{\Gamma(k+n+2\lambda)}.
\end{split}
\end{equation}
By fundamental theorem of gamma function, we see that
\begin{equation}\label{43}
\Gamma\left(k+\lambda+l+\frac{1}{2}\right)=\binom{k+\lambda+l-\frac{1}{2}}{l}l!\Gamma\left(k+\lambda+\frac{1}{2}\right),
\end{equation}
\begin{equation}\label{44}
\Gamma\left(\lambda+n-l+\frac{1}{2}\right)=\binom{\lambda+n-l-\frac{1}{2}}{n-l}(n-l)!\Gamma\left(\lambda+\frac{1}{2}\right),
\end{equation}
and
\begin{equation}\label{45}
\Gamma\left(k+2\lambda+n\right)=\binom{k+2\lambda+n-1}{n+k}(n+k)!\Gamma(2\lambda).
\end{equation}
As is well known, the duplication formula for  gamma function is given by
\begin{equation}\label{46}
\Gamma(z)\Gamma\left(z+\frac{1}{2}\right)=2^{1-2z}\sqrt{\pi}\Gamma(2z).
\end{equation}
By \eqref{42}, \eqref{43}, \eqref{44} and \eqref{45}, we get
\begin{equation}\label{47}
\begin{split}
&\int_{-1} ^1(1-x)^{k+\lambda+l-\frac{1}{2}}(1+x)^{\lambda+n-l-\frac{1}{2}}dx \\
=&2^{k+n+1}\frac{\binom{k+\lambda+l-\frac{1}{2}}{l}\binom{\lambda+n-l-\frac{1}{2}}{n-l}\Gamma\left(k+\lambda+\frac{1}{2}\right)}{\binom{n}{l}\binom{k+2\lambda+n-1}{n+k}\binom{n+k}{k}k!\Gamma(\lambda)}\sqrt{\pi}.
\end{split}
\end{equation}
From \eqref{41} and \eqref{47}, we have
\begin{equation}\label{48}
\begin{split}
d_k&=\lambda^k(k+\lambda)2^{2k+1}\frac{\binom{n+k+2\lambda-1}{n-k}}{\binom{n+\lambda-\frac{1}{2}}{n-k}} \\
&\times \sum_{l=0} ^{n-k}\binom{n+\lambda-\frac{1}{2}}{n-k-l}\binom{n+\lambda-\frac{1}{2}}{l}(-1)^l\frac{\binom{k+\lambda+l-\frac{1}{2}}{l}\binom{\lambda+n-l-\frac{1}{2}}{n-l}}{\binom{n}{l}\binom{k+2\lambda+n-1}{n+k}\binom{n+k}{k}k!}.
\end{split}
\end{equation}
Therefore, by \eqref{48}, we obtain the following theorem.
\begin{theorem}
For $n \in {\mathbb{Z}}_+$, we have
\begin{equation*}
\begin{split}
C_n ^{(\lambda)}(x)&=\sum_{k=0} ^n \left\{\frac{\lambda^k(k+\lambda)2^{2k+1}\binom{n+k+2\lambda-1}{n-k}}{\binom{n+\lambda-\frac{1}{2}}{n-k}} \right.\\
 & \left. \times \sum_{l=0} ^{n-k}
 \frac{\binom{n+\lambda-\frac{1}{2}}{n-k-l}\binom{n+\lambda-\frac{1}{2}}{l}(-1)^l \binom{k+\lambda+l-\frac{1}{2}}{l} }{\binom{n}{l}\binom{k+2\lambda+n-1}{n+k}\binom{n+k}{k}k!}\binom{\lambda+n-l-\frac{1}{2}}{n-l}\right\}C_k ^{(\lambda)}(x).
\end{split}
\end{equation*}
\end{theorem}

\end{document}